\documentclass[letterpaper, paper,11pt]{AAS}	% for preprint proceedings
\usepackage{amsfonts}
\usepackage{bm}
\usepackage{amsmath}
\usepackage{subfigure}
\usepackage[colorlinks=true, pdfstartview=FitV, linkcolor=black, citecolor= black, urlcolor= black]{hyperref}
\usepackage{footnpag} % make footnote symbols restart on each page	 
\usepackage{subcaption,multicol,lipsum}
\PaperNumber{24-302}

\newcommand{\norm}[1]{\left\lVert#1\right\rVert}

\begin{document}

%\title{A comparison between L2-norm-based and Hyperbolic-Tangent-Based Smoothing With State Transition Matrix for Low-thrust trajectories }

\title{Comparison of control regularization techniques for minimum-fuel low-thrust trajectory design using indirect methods}

\author{Saeid Tafazzol\thanks{Graduate Student, Department of Aerospace Engineering, Auburn University, 141 Engineering Dr Auburn, AL 36849.}, Ehsan Taheri\thanks{Assistant Professor, Department of Aerospace Engineering, Auburn University, 141 Engineering Dr Auburn, AL 36849.}, Nan Li\thanks{Assistant Professor, Department of Aerospace Engineering, Auburn University, 141 Engineering Dr Auburn, AL 36849.}}

\maketitle{} 	

\begin{abstract}
Minimum-fuel low-thrust trajectories typically consist of a finite, yet unknown number of switches in the thrust magnitude profile. This optimality-driven characteristic of minimum-fuel trajectories poses a challenge to the numerical methods that are typically used for solving the resulting Hamiltonian boundary-value problems (BVPs). In this paper, we compare the impact of the popular hyperbolic-tangent-based smoothing with a novel L2-norm-based smoothing on the convergence performance of quasi-Newton gradient-based methods. Both smoothing methods are applied directly at the level of control, which offer a significant implementation simplicity. We also investigate the application of each method in scenarios where the State Transition Matrix (STM) is used for accurate calculation of the sensitivities of the residual vector of the resulting BVPs with respect to the unknown initial costate values. The effectiveness of each control smoothing method is assessed across several benchmark minimum-fuel low-thrust trajectory optimization problems.

%In many optimal control problems (OCP), the solution is known as ``bang-bang'' indicating that the optimal control alternates between discrete values at specific points in time. This characteristic .  In this paper, we compare the performance of Hyperbolic-tangent-based smoothing to that of a recently introduced L2-norm-based smoothing. Additionally, we examine the application of each method in scenarios where the State Transition Matrix (STM) is employed. The capability of each method to solve optimal problems is demonstrated through the example of low-thrust fuel-optimal trajectories.
\end{abstract}

\section{Introduction}

Trajectory optimization of space vehicles \cite{whiffen2006mystic} is an important task in space mission design and for solving problems in orbital mechanics and astrodynamics  \cite{ellison2018application}. Efficient techniques that are used for solving trajectory optimization problems affect the quality and economic cost of space missions. In the past couple of decades, advances in electric propulsion technology have led to the design and use of spacecraft thrusters. Due to their high specific impulse values (typically one order of magnitude larger than their chemical rocket counterparts), electric thrusters offer significant fuel efficiency  \cite{lev2019technological}. 

Upon using electric propulsion systems, the resulting low-thrust trajectories may consist of multiple orbital revolutions around the central body \cite{petukhov2019application,meng2019low}. In fact, for low-thrust planet-centric trajectories, the number of orbital revolutions can exceed beyond a few hundreds \cite{aziz2018low,taheri2021optimization}. In addition, depending on the considered performance index (i.e., when optimality is sought with respect to either fuel consumption or time of flight), the spacecraft thruster(s) may not operate at all times and the optimal solution may consist of a finite, but unknown number of thrusting and coasting arcs. The sequence of thrusting and coasting arcs is referred to as the structure of the optimal control profile, where the duration of each one of the arcs is unknown and has to be determined. 

Numerical methods that are used for solving practical optimal control and trajectory optimization problems are broadly classified under \textit{direct} and \textit{indirect} optimization methods \cite{betts1998survey,trelat2012optimal}. Both numerical methods have been used extensively for solving spacecraft trajectory optimization problems. Indirect/variational approaches to optimal control, rooted in the calculus of variations, construct a boundary-value problem (BVP) by using Lagrangian multipliers (also called co-states) associated with every state to incorporate the system's dynamics as nonlinear equality constraints, resulting in Hamiltonian BVPs. The simplest form of BVPs, i.e., when no control or state-path constraints are considered, reduce to two-point boundary value problems (TPBVPs). To solve the resulting TPBVPs, one can leverage quasi-Newton gradient-based root-solving methods that iterate (for most problems) on the unknown co-state values at boundaries to generate solutions (consisting of states, costates, and control time histories) that satisfy necessary conditions of optimality. In this paper, we focus on indirect optimization methods, for three reasons: 
\begin{itemize}
    \item 1) In astrodynamics and orbital mechanics problems, the motion of spacecraft is predominantly governed by the Keplerian motion (for the most part of the trajectory within the sphere of influence of the central body). Thus, the acceleration produced by low-thrust propulsion systems is small compared to the dominant central-body gravitational acceleration and third-body perturbations. Even if the control produced by the propulsion system is considered, the changes of the orbital elements remain relatively small (depending on the magnitude of the thrust). This makes low-thrust trajectories trajectory optimization problems to become ``easier'' to solve, compared to the trajectory optimization of atmospheric flight vehicles that are subject to significant forces that can no longer be considered as perturbations. Thus, indirect methods remain an alternative choice for formulating and solving low-thrust trajectory optimization problems. 
    \item 2) Indirect approaches, in principle, produce high-resolution solutions that represent the most accurate achievable solutions for the same local/global extremal/optimal solutions. In the context of trajectory optimization, high-resolution solutions correspond to ones that present changes in the states and controls with a high temporal resolution, and more importantly, capture all events (such as entry to and exit from eclipses) with a high accuracy. This accuracy is important since one can use the solutions obtained using indirect methods as baseline solutions to determine the gap (or sacrifice) in optimality (however small or large) associated with the solutions obtained using direct optimization methods \cite{ross2007low}. An example comparison between indirect and convex-optimization based methods is conducted in \cite{nurre2022comparison} when minimum-fuel trajectories are considered. Perhaps, another equally important distinction is that historically, indirect optimization methods have been combined with homotopy and numerical continuation methods. This combination has been a necessity to alleviate the numerical challenges associated with solving practical optimal control problems, but homotopy techniques have revealed a more ``global'' and insightful mapping of the space of solutions (a.k.a. extremal field maps) \cite{taheri2020many}. Examples include the class of minimum-time orbit transfer maneuvers with many local extremals \cite{pan2019finding,pan2020practical,pan2020bounding,zhang2023solution}. For the class of fixed-time rendezvous-type minimum-fuel trajectory optimization problems, a systematic study can be performed to the number of orbital revolutions to determine the global minimum-fuel solution \cite{taheri2020many}. 
    
    \item 3) The idea of control regularization, which falls under the general principle of ``invariant embedding'' is a key enabler for solving TPBVPs. In particular, consideration of state, control and mixed state-control inequality constraints is more relevant and of interest, when we use indirect formalism of optimal control theory. Direct optimization methods, have been historically favored since they can deal with different types of constraints due to the fact that constraints are enforced at a finite number of (mesh) points and have a larger domain of convergence compared to indirect methods \cite{malyuta2021advances}. For indirect methods, numerical solutions of practical optimal control problems have been faced with various challenges associated with not only discontinuous changes in controls, but also due to derivation of the necessary (most of the time) conditions of optimality, in particular, for problems that have to deal with different types of inequality constraints. However, significant progresses have been made in recent years that have allowed indirect methods to solve quite challenging optimal control problems including the unified trigonometrizaton method, UTM \cite{taheri2020minimum,mall2020entry,mall2020unified,mall2020uniform,mall2022three}, vectorized trigonometric regularization (VTR) \cite{kovryzhenko2023vectorized}, and generalized vectorized trigonometric regularization (GVTR) \cite{kovryzhenko2024generalized}. %The paper is related to the recent progresses made in control regularization methods. 
\end{itemize}

In minimum-fuel trajectory optimization problems, the presence of thrusting and coasting arcs is a consequence of applying Pontryagin's minimum principle  \cite{haberkorn2004low}. For minimum-time maneuvers, the optimality principle requires the thruster to operate at its maximum thrusting capability during the entire maneuver time \cite{taheri2017co}; however, it is still possible to have coasting arcs if the optimization problem has to take into account eclipses mostly during the planet-centric phases of flights \cite{cerf2019fast,sowell2024eclipse} and for cislunar trajectory optimization problems \cite{nurre2024end}. For minimum-fuel trajectory optimization problems, the optimal bang-off-bang thrust profile causes the Jacobian (used in gradient-based solvers for determining the search direction) to become singular \cite{bertrand2002new,haberkorn2004low}, which in turn deteriorates the selection of the descent direction and renders the numerical solution impossible. %The piece-wise continuous optimal throttle structure is a source of difficulty when numerical methods are used, since the Jacobian of the residual vector with respect to the unknown decision values becomes singular \cite{}. 
The existence of multiple throttle switches exacerbates the issue, and reduces the basin of attraction of the numerical methods \cite{haberkorn2004low,prussing2010primer,jiang2012practical}. 

To overcome the singularity in the Jacobian, the TPBVPs associated with minimum-fuel trajectory optimization problems are typically solved through various regularization methods, such as those that regularize the non-smooth thrust/throttle profiles. Control regularization methods such as logarithmic smoothing \cite{bertrand2002new}, extended logarithmic smoothing \cite{taheri2016enhanced}, hyperbolic tangent smoothing (HTS) \cite{taheri2018generic} are among the popular and efficient methods for attenuating the issues with non-smooth control profiles. Through regularization methods, the non-smooth optimal control problem is embedded into a one- or multiple-parameter family of smooth neighboring optimal control problems. Then, standard (or advanced) numerical continuation and/or homotopy methods are used to solve the resulting family of smooth TPBVPs until a solution to the original optimal control problem is obtained \cite{trelat2012optimal}. 

In terms of applications, the HTS method is already used to solve interplanetary minimum-fuel trajectory optimization and minimum-time rest-to-rest satellite reorientation problems \cite{taheri2018generic}. One of the main advantages of the HTS method is its ease of implementation. One can proceed with the standard approach to the optimal control problem and determine the so-called thrust/throttle switching function. Then, the HTS method is applied as a filter to the thrust switching function to form a smooth approximation of the theoretically-optimal bang-bang thrust profile. Thus, the HTS method is directly applied at the control level, which is an important feature because the standard procedure for formulating the necessary conditions (using indirect methods) remains unaffected. The HTS method is also a key component of a novel framework -- Composite Smooth Control (CSC) -- proposed for solving optimal control problems with discrete and multiple modes of operations \cite{taheri2020novel,TAHERI2020166}. The CSC framework and the HTS method are also used for co-optimization of spacecraft propulsion system and its trajectory \cite{arya2021composite}, optimal mode-selection and net payload mass optimization using indirect optimization methods \cite{arya2021electric}, multimode trajectory optimization with low-thrust gravity-assist maneuvers \cite{arya2021low}, eclipse-conscious cislunar low-thrust trajectory design \cite{singh2021eclipse}, and for minimum-fuel asteroid landing trajectory optimization problems \cite{nakano2022time}. 

In this paper, we compare two state-of-the-art regularization methods, i.e., the HTS method and a recent L2-norm-based regularization, proposed by Taheri and Li \cite{taheri2023l2}. In \cite{taheri2023l2}, the utility of the L2-norm-based regularization is demonstrated for solving three optimal control problems: 1) minimum-fuel low-thrust trajectories, 2) the standard Goddard rocket problem with its characteristic bang-singular-bang thrust profile and 3) a minimum-time spacecraft reorientation problem with both bang-bang
and a second-order singular arc. The L2-norm-based regularization is already used for co-optimization of spacecraft propulsion system and its trajectory using a direct optimization method \cite{saloglu2024co}. The HTS method is already compared against the logarithmic smoothing in \cite{taheri2018generic} and against the standard quadratic smoothing in \cite{taheri2018performance} revealing its advantages by broadening the basin of attraction of the resulting TPBVPs. Application of the HTS combined with the State Transition Matrix (STM) is also studied in \cite{arya2019hyperbolic}. Similar to \cite{arya2019hyperbolic}, we consider the impact of accurate calculation of sensitivities for the HTS and L2-norm-based regularization/smoothing methods, when finite-difference and the STM methods, are used. We compare the results of the two methods on two benchmark minimum-fuel trajectory optimization problems. 

\section{Problem Formulation}
To assess the impact of control regularization on the convergence performance of numerical methods, we consider two different parameterizations of the motion dynamics using Cartesian coordinates and the set of modified equinoctial elements \cite{taheri2016enhanced,junkins2019exploration}. In addition, we consider optimizing the trajectory of spacecraft during heliocentric phases of flight with zero hyperbolic excess velocity relative to the departing and arrival bodies. We consider two-body dynamics assumptions and ignore third-body and other types of perturbations (e.g., solar-radiation pressure). The only perturbing acceleration (and at the same time control input) is due to the operation of the propulsion system during the thrusting arcs. These assumptions are consistent with the dynamical models and boundary conditions of the chosen benchmark problems \cite{taheri2016enhanced}. 
\subsection{Equations of motion}
Following \cite{junkins2019exploration}, we adopt a unified control-affine representation of the dynamics. Let $\bm{x} \in \mathbb{R}^6$ denote part of the state vector (e.g., Cartesian coordinates consisting of components of the position and velocity vectors or the six modified equinoctial elements). Let $m \in \mathbb{R}^+$ denote the instantaneous mass of the spacecraft. We adopt a magnitude-vector parameterization of the control-acceleration vector. Let $\delta \in [0,1]$ denote  engine throttle input and let $\hat{\bm{\alpha}}$ denote the thrust steering unit vector (i.e.,  $\| \hat{\bm{\alpha}}\| = 1$), the acceleration vector produced by the propulsion system can be written as,
\begin{align}
\bm{u} = \frac{T_\text{max}}{c} \hat{\bm{\alpha}}\delta, 
\end{align}
where $T_\text{max}$ denotes the maximum thrust magnitude and $c = I_\text{sp} g_0$ denotes the effective exhaust velocity. Here, $I_\text{sp}$ is the thruster specific impulse value and $g_0$ is the Earth's sea-level gravitational constant. To simplify the problem formulation, it is assumed that the values of $T_\text{max}$ and $c$ remain constant throughout the trajectory. In reality, and for solar-powered low-thrust propulsion systems, these parameters depend on the power that is sent to the engine \cite{arya2021composite}. Altogether, the control vector can be written as $\bm{u}^\top = [\hat{\bm{\alpha}}^\top,\delta] \in \mathbb{R}^4$. The time rate of change of $\bm{x}$ and $m$ can be written as,
\begin{align} \label{eq:statemassdynamics}
    \dot{\bm{x}}(t) = \mathbb{A}(\bm{x}(t),t) + \mathbb{B}(\bm{x}(t),t) \left (\frac{T_\text{max}}{m} \hat{\bm{\alpha}}\delta \right ), \quad \dot{m} = -\frac{T_\text{max}}{c} \delta,
\end{align}
where the entries of the $\mathbb{A}$ vector and the $\mathbb{B}$ matrix depend on the choice of coordinates and/or elements \cite{taheri2016enhanced,junkins2019exploration}. The main task, in spacecraft trajectory optimization problems, is to determine the time history of the optimal/extremal (denoted by superscript `*') control vector, $\bm{u}^*$, such that a performance index is minimized/maximized while equations of motions, given in Eq.~\eqref{eq:statemassdynamics}, are satisfied along with additional boundary conditions and along-the-path state and control equality/inequality constraints. We proceed by defining the cost functional (a.k.a. the performance index or objective) of the minimum-fuel trajectory optimization problems.
\subsection{Cost functional of minimum-fuel trajectory optimization problems}
In this paper, we focus on minimum-fuel trajectory optimization problems. In optimal control terminology, the cost functional, $J \in \mathbb{R}$, can be stated in various equivalent forms \cite{bryson2018applied}. However, in this paper, we consider the Lagrange form in which the cost functional is expressed in the form of an integral. For spacecraft with a fixed initial mass, minimum-fuel trajectories correspond to minimizing the propellant consumption, which can be stated mathematically as,
\begin{align}
    \underset{\delta \in [0,1]~\&~ \hat{\bm{\alpha}} \in \mathcal{U}_{\hat{\bm{\alpha}}}}{\text{Minimize}} \quad J = \int_{t_0}^{t_f} \frac{T_\text{max}}{c} \delta(t)~dt,
\end{align}
where $\mathcal{U}_{\hat{\bm{\alpha}}} = \{\hat{\bm{\alpha}}: [t_0,+\infty) \rightarrow \bm{U} ~|~ \hat{\bm{\alpha}}~\text{is measurable} \}$ with $\bm{U} = \{ \hat{\bm{\alpha}}~|~\| \hat{\bm{\alpha}} \| = 1 \}$. The admissible set of throttle control is listed under the minimization operator. 

Please note that while optimization appears to be performed only over the control vector, $\bm{u}$, the optimization is actually performed over all admissible sets of states and controls. However, the time derivative of states (consisting of $\bm{x}$ and $m$) is governed by the right-hand side of the differential equations that is given in Eq.~\eqref{eq:statemassdynamics}. Please also note that one could have formulated the cost functional in Mayer form, as $J = -m(t_f)$ that considers the negative of the final mass of the spacecraft, minimization of which is equivalent to maximizing the final mass. We proceed by writing the optimal control problem formulation associated with the minimum-fuel trajectory optimization problems, when Cartesian coordinates and MEEs are used. 

%The necessary condition for the fuel-optimal trajectory can be obtained equivalently by maximizing the final mass, i.e., $m(t_f)$. By incorporating Eq.~\eqref{eom_cart}c, we can obtain the following functional:

%\begin{equation}
%\begin{aligned}
%&\text{Minimize}\; J = \dfrac{T_{\max}}{c} %\int_{t_0}^{t_f} \delta \, dt \\
%&\text{S.t. }\\
%&\mathbf{r}(t_0) - \mathbf{r}_0 = \mathbf{0},\\
%&\mathbf{v}(t_0) - \mathbf{v}_0 = \mathbf{0},\\
%&\mathbf{r}(t_f) - \mathbf{r}_f = \mathbf{0},\\
%&\mathbf{v}(t_f) - \mathbf{v}_f = \mathbf{0},\\
%&\dot{\mathbf{x}} = \bm{f}(\mathbf{x},\bm{\alpha},t),\\
%\end{aligned}
%\end{equation}
\subsection{Formulation of the minimum-fuel problem using Cartesian coordinates}
Let $\bm{r} = [x,y,z]^\top$ and $\bm{v} = [v_x,v_y,v_z]^\top$ denote the position and velocity vectors of the center of mass of the spacecraft relative to the origin of the heliocentric frame of reference. Thus, we have $\bm{x}^\top = [\bm{r}^\top,\bm{v}^\top]$.  The time rate of change of states can be written as, \cite{junkins2019exploration}
\begin{align}\label{eom_cart}
\dot{\bm{r}} &= \bm{v} = \bm{f}_{\bm{r}}, \nonumber \\
\dot{\bm{v}} &= -\dfrac{\mu}{r^3} \bm{r} + \dfrac{T_{\max}}{m} \hat{\bm{\alpha}}  \delta= \bm{f}_{\bm{v}}, \nonumber \\
\dot{m} &= - \dfrac{T_{\max}}{c} \delta = f_m,  
\end{align}
where $\mu$ is the gravitational parameter for the central body (Sun in our problems) and $r = \|\bm{r} \|$. The time rate of change of position and velocity vectors, given in Eq.~\eqref{eom_cart}, can be written as $\dot{\bm{x}} = \bm{f}(\bm{x}(t),\bm{u}(t),t)$. The minimum-fuel optimal control problem can be written as,
\begin{equation}
\begin{aligned}
&\underset{\delta \in [0,1]~\&~ \hat{\bm{\alpha}} \in \mathcal{U}_{\hat{\bm{\alpha}}}}{\text{Minimize}}\; J =  \int_{t_0}^{t_f} \dfrac{T_{\max}}{c}\delta(t) \, dt \\
&\text{s.t.:}\\
&\dot{\bm{x}}(t) = \bm{f}(\bm{x}(t),\bm{u}(t),t), \quad \dot{m} = -\frac{T_\text{max}}{c} \delta,\\
&\bm{r}(t_0) - \bm{r}_0 = \bm{0},\\
&\bm{v}(t_0) - \bm{v}_0 = \bm{0},\\
&\bm{r}(t_f) - \bm{r}_f = \bm{0},\\
&\bm{v}(t_f) - \bm{v}_f = \bm{0},\\
& m(t_0) = m_0,
\end{aligned}
\end{equation}
where $\bm{x}_0^\top =[\bm{r}_0^\top, \bm{v}_0^\top]$ and $m_0$ denote the complete states of the spacecraft at the time of departure from the departing planet (Earth) and $\bm{r}_f$ and $\bm{v}_f$ denote the fixed position and velocity vectors of the target/arrival body. 
We can use the indirect formalism of optimal control theory, which treats the state differential equation as constraints \cite{bryson2018applied}. Let $\bm{\lambda} ^\top= [\bm{\lambda}_{\bm{r}}^\top, \bm{\lambda}_{\bm{v}}^\top]$ denote the vector of Cartesian coordinates and let $\lambda_{m}$ denote the mass costate. We have dropped the argument $t$ from the terms to make the relations more concise. We can proceed by forming the variational Hamiltonian as, 
\begin{equation}
H = \dfrac{T_{\max}}{c} \delta + \boldsymbol{\lambda}_{\bm{r}}^\top \bm{f}_{\bm{r}} + \boldsymbol{\lambda}_{\bm{v}}^\top \bm{f}_{\bm{v}} + \lambda_m f_m
\end{equation}

The differential equations for the costates are determined using the Euler-Lagrange equation \cite{bryson2018applied} as,
\begin{align} \label{eq:ELcar}
    \dot{\bm{\lambda}} & = -\left[\dfrac{\partial H}{\partial \bm{x}}\right]^\top, & \dot{\lambda}_m & = -\frac{\partial H}{\partial m}.
\end{align}

Following the primer vector theory of Lawden \cite{prussing2010primer} and Pontryagin's minimum principle, the extremal control expressions \cite{taheri2018performance,arya2019hyperbolic} can be derived and stated compactly as, 
\begin{align} \label{eq:optimalcar}
    \hat{\bm{\alpha}}^* & = -\frac{\bm{\lambda_{\bm{v}}}}{\|\bm{\lambda_{\bm{v}}}\|}, & \delta^* & \begin{cases}
    =  1      & \text{if } S > 0, \\
    \in   [0,1] & \text{if } S = 0, \\
    =  0      & \text{if } S \leq 0,
\end{cases} & \text{with} \quad  S = \frac{c\|\bm{\lambda}_{\bm{v}}\|}{m} + \lambda_m - 1,
\end{align}
where $S$ is the so-called throttle \textit{switching function}. The primer vector is $\bm{p} = -\bm{\lambda_v}$.
%\begin{equation}
%    \delta^* = \arg \min_{\delta \in [0,1]} H(\boldsymbol{x}^*(t), \delta(t), \boldsymbol{\lambda}^* (t))
%\end{equation}
%To solve for \(\delta\), the Hamiltonian can be written as:
%\begin{subequations}
%\begin{align}
%H &= H_0 - H_\delta \delta  \\
%H_0 &= \bm{\lambda}_{\bm{r}}^T v + \bm{\lambda}_{\bm{v}}^T \left( -\frac{\mu}{r^3} \bm{r} \right),  \\
%H_\delta &= \frac{T_{\text{max}}}{c} \left[ \lambda_m + \frac{c}{m} \bm{\lambda}_{\bm{v}}^T \bm{\alpha} - 1 \right]. 
%\end{align}
%\end{subequations}
%Where $H_0$ is free of $\delta$. Therefore, by solving this problem using PMP \cite{taheri2018performance,arya2019hyperbolic}, we obtain the following expression for $\delta^*$:
%\begin{equation}
%    SF = \frac{\|\boldsymbol{\lambda}_{\boldsymbol{v}}\| c}{m} + \lambda_m - 1
%\end{equation}
%where $SF$ is the so-called thrust switching function. 
According to Lawden, and for problems in which the direction of the thrust vector is not constrained, the time history of $\bm{\lambda}_{\bm{v}}$ is continuous, which results in a continuous profile of $\hat{\bm{\alpha}}^*$. However, the value of $\delta^*$ depends on the sign of the throttle switching function, which results in discontinuities in $\delta^*$. Note that the number and time instants of the zero-crossings of $S$ are not known \textit{a priori}. The presence of instantaneous changes in the value of $\delta^*$ poses significant challenges to numerical solvers and integrators. To overcome this issue, smoothing functions are used, which employ squeezing continuous functions such as the hyperbolic tangent function to map the switching function between 0 and 1. In this paper, we compare the following smoothing functions:
\begin{enumerate}
    \item Hyperbolic-Tangent-Based Smoothing \cite{taheri2018generic}:
    \begin{equation}
        \delta_\text{tanh}^* \approx \delta_\text{tanh} (S;\rho) =  0.5 \left[1 + \tanh\left(\frac{S}{\rho}\right) \right],
        \label{eqn:smooth_tanh}
    \end{equation}
    
    \item L2-norm-Based Smoothing \cite{taheri2023l2}:
    \begin{equation}
        \delta^*_\text{L2} \approx \delta_\text{L2} (S;\rho) = 0.5 \left[1 + \frac{S}{\sqrt{S^2 + \rho^2}} \right],
        \label{eqn:smooth_l2}
    \end{equation}
\end{enumerate}
where $\rho$ denotes a smoothing parameter that allows us to control the sharpness of change in throttle. This parameter will be used within a standard numerical continuation method. The main idea, in control regularization, is to embed the non-smooth, piece-wise continuous optimal throttle logic, given in Eq.~\eqref{eq:optimalcar}, into a one-parameter family of smooth curves. In the limit and as the value of $\rho$ is decreased (from a relatively large value, e.g., $\rho = 1$) to small values (e.g., $\rho = 1.0 \times 10^{-3}$), the smooth throttle approaches the piece-wise continuous throttle, given in Eq.~\eqref{eq:optimalcar}. 

The resulting one-parameter family of smooth TPBVPs is summarized as follows: the state differential equations, given in Eq.~\eqref{eom_cart}, the costate differential equations, obtained from Eq.~\eqref{eq:ELcar}, the extremal thrust steering unit vector, $\hat{\bm{\alpha}}$ and the throttle switching function, given in Eq.~\eqref{eq:optimalcar}, and the smooth throttle function (either HTS or L2-norm based), given in Eq.~\eqref{eqn:smooth_tanh} or \eqref{eqn:smooth_l2}, and the transversality condition on the costate associated with mass at the final time, $\lambda_m(t_f) = 0$. 

The resulting TPBVP is written as a nonlinear shooting problem as,
\begin{align}
    \bm{\psi}(\bm{\eta}(t_0);\rho) = [\bm{r}^\top(t_f)-\bm{r}^\top_f,\bm{v}^\top(t_f)-\bm{v}^\top_f,\lambda_m(t_f)] = \bm{0}_{7\times 1},
\end{align}
where $\bm{\eta}^\top(t_0) = [\bm{\lambda_{r}}^\top(t_0),\bm{\lambda_{v}}^\top(t_0),\lambda_m(t_0)] \in \mathbb{R}^7$ denotes the vector of unknown initial costates. Typically, the resulting shooting problems are solved using single- or multiple-shooting solution schemes \cite{nurre2024end}. Note that the initial states at $t = t_0$ are known and these conditions are not added to the boundary conditions of the shooting problem. Numerical continuation methods are typically used for solving the problems by setting $\rho = 1$ and solving the shooting problem. Once a solution to $\bm{\eta}^\top(t_0)$ is obtained, the value of $\rho$ is decreased, say by a factor of 0.1 or 0.5 (i.e., $\rho = \rho \times 0.1$). The previous value of $\bm{\eta}^\top(t_0)$ is used as an initial guess for the slightly-modified shooting problem. These steps are repeated until the value of $\rho$ is smaller than some prescribed user-defined value or subsequent changes in the value of cost, $J$, becomes smaller than a threshold.

%\subsection{State Transition Matrix}
\subsection{Formulation of the minimum-fuel problem using modified equinoctial elements}
Let the state vector associated with the MEEs be denoted as $\bm{x} = [p,~f,~g,~h,~k,~L]^\top$ (with an abuse of notation). The mappings between MEEs and the Keplerian classical orbital elements (COEs) are as follows \cite{walker1985set}: $p = a(1 - e^2)$,
    $f = e \cos{(\omega + \Omega)}$, 
    $g = e \sin{(\omega + \Omega)}$,
    $h = \tan\left(\tfrac{i}{2}\right)\cos(\Omega)$,
    $k = \tan\left(\tfrac{i}{2}\right)\sin(\Omega)$, $
    L = \theta + \omega + \Omega$, where $a$, $e$, $i$, $\Omega$, $\omega$ and $\theta$, denote the semi-major axis, eccentricity, inclination, right-ascension of the ascending note, argument of periapses, and true anomaly, respectively. The MEEs present an ideal parameterization of the evolution of low-thrust trajectories, avoiding the singularities (associated with zero inclination and eccentricity values) present in the classical orbital elements \cite{junkins2019exploration}. The time rate of change of MEEs and mass can be written as,
\begin{align}  \label{eom_MEE}  
    \dot{\bm{x}} &= \mathbb{A} +\mathbb{B} ( \frac{T_{\max}}{m} \hat{\bm{\alpha}} \delta ) & \dot{m} & = -\frac{T_\text{max}}{c} \delta,
    \end{align}
with $\mathbb{A}$ and $\mathbb{B}$ defined (under the canonical unit scaling with $\mu = 1$) as,
\begin{align}  
\mathbb{A}^\top &= \begin{bmatrix}
        0 &
        0 &
        0 &
        0 &
        0 &
        \sqrt{p} \left( \frac{q}{p} \right)^2 
    \end{bmatrix}, \\
     \mathbb{B} &= \begin{bmatrix}
        0 & \frac{2p}{q} \sqrt{p} & 0 \\
        \sqrt{p} \sin L & \frac{\sqrt{p}}{q} \left( (q + 1) \cos L + f \right) & -\frac{\sqrt{p} g}{q} \left( h \sin L - k \cos L \right) \\
        -\sqrt{p} \cos L & \frac{\sqrt{p}}{q} \left( (q + 1) \sin L + g \right) & \frac{\sqrt{p} f}{q} \left( h \sin L - k \cos L \right) \\
        0 & 0 & \frac{\sqrt{p} s^2 \cos L}{2q} \\
        0 & 0 & \frac{\sqrt{p} s^2 \sin L}{2q} \\
        0 & 0 & \frac{\sqrt{p}}{q} \left( h \sin L - k \cos L \right) \\
    \end{bmatrix}, 
   \end{align}
with $s^2 = 1 + h^2 + k^2$ and $q = 1 + f \cos L + g \sin L$.

The minimum-fuel optimal control problem can be written as,
\begin{equation}
\begin{aligned}
&\underset{\delta \in [0,1]~\&~ \hat{\bm{\alpha}} \in \mathcal{U}_{\hat{\bm{\alpha}}}}{\text{Minimize}}\; J =  \int_{t_0}^{t_f} \dfrac{T_{\max}}{c}\delta(t) \, dt \\
&\text{s.t.:}\\
&\dot{\bm{x}}(t) = \mathbb{A}+\mathbb{B} (\frac{T_{\max}}{m} \hat{\bm{\alpha}} \delta),\quad \dot{m}(t) = - \dfrac{T_{\max}}{c} \delta ,\\
&\bm{x}(t_0) - \bm{x}_0 = \bm{0},\\
&\bm{x}(t_f) - \bm{x}_f = \bm{0},\\
& m(t_0) = m_0,
\end{aligned}
\end{equation}
where $\bm{x}_0^\top =[p_0, f_0, g_0, h_0, k_0, L_0]$ denotes the MEEs of the spacecraft at the time of departure from the departing planet (Earth) and $\bm{x}_f^\top =[p_f, f_f, g_f, h_f, k_f, L_f]$ denote the fixed vector of target MEE values at the end of the trajectory. 

% \subsection{Cost Functional}
Let $\bm{\lambda}$ denote (with an abuse of notation) the costate vector associated with the MEEs and let $\lambda_m$ denote the costate associated with mass. We can proceed and form the Hamiltonian as,
\begin{align}
    H = \frac{T_\text{max}}{c} \delta + \bm{\lambda}^{\top} (\mathbb{A} + \mathbb{B} (\frac{T_{\max}}{m} \delta\bm{\alpha} ) ) - \lambda_m \frac{T_\text{max}}{c} \delta.
\end{align}

The differential equations for the costates associated with MEEs are determined using the Euler-Lagrange equation \cite{bryson2018applied} as,
\begin{align} \label{eq:ELMEE}
    \dot{\bm{\lambda}} & = -\left[\dfrac{\partial H}{\partial \bm{x}}\right]^\top, & \dot{\lambda}_m & = -\frac{\partial H}{\partial m}.
\end{align}

Following the steps that we followed, when we used the Cartesian coordinates, and upon using the primer vector theory of Lawden and the PMP, the extremal control expressions can be derived and written as
\begin{align} \label{eq:optimalMEE}
    \hat{\bm{\alpha}}^* & = -\frac{\mathbb{B}^\top \bm{\lambda}}{\|\mathbb{B}^\top \bm{\lambda}\|}, & \delta^* & \begin{cases}
    =  1      & \text{if } S > 0, \\
    \in   [0,1] & \text{if } S = 0, \\
    =  0      & \text{if } S \leq 0,
\end{cases} & \text{with} \quad  S = \frac{c\|\mathbb{B}^\top \bm{\lambda}\|}{m} + \lambda_m - 1,
\end{align}
where $S$ is the so-called throttle \textit{switching function}. It is important to notice that the primer vector, is obtained by modulating the costate vector associated with MEEs by the control influence matrix, i.e., $\bm{p} = -\mathbb{B}^\top \bm{\lambda}$. Please also note that the primer vector, when Cartesian coordiates are used, can be expressed in a similar form, i.e., $\bm{p} = -\mathbb{B}^\top \bm{\lambda}$ except that the control influence matrix for the set of Cartesian coordinates is $\mathbb{B} = [\bm{0}_{3\times3}, I_{3\times3}]$ and the costate vector is $\bm{\lambda} = [\bm{\lambda_r}^\top,\bm{\lambda_v}^\top]$. Thus, we have $\bm{p} = -\mathbb{B}^\top \bm{\lambda} = -\bm{\lambda_v}$, which is identical to the primer vector relation given below Eq.~\eqref{eq:optimalcar}. 

The resulting smooth TPBVP is summarized as follows: the state differential equations, given in Eq.~\eqref{eom_MEE}, the costate differential equations, obtained from Eq.~\eqref{eq:ELMEE}, the extremal thrust steering unit vector, $\hat{\bm{\alpha}}$ and the throttle switching function, given in Eq.~\eqref{eq:optimalMEE}, and the smooth throttle function (either HTS or L2-norm based), follow the same relations given in Eq.~\eqref{eqn:smooth_tanh} or \eqref{eqn:smooth_l2} except that the switching function is calculated using the relation given in Eq.~\eqref{eq:optimalMEE}, and the transversality condition on the costate associated with mass at the final time, $\lambda_m(t_f) = 0$. 

The resulting TPBVP is written as a shooting problem as
\begin{align}
    \bm{\psi}(\bm{\eta}(t_0);\rho) = [\bm{x}^\top(t_f)-\bm{x}_f^\top,\lambda_m(t_f)] = \bm{0}_{7\times 1},
\end{align}
where $\bm{\eta}^\top(t_0) = [\bm{\lambda}^\top(t_0),\lambda_m(t_0)] \in \mathbb{R}^7$ denotes the vector of unknown initial costates. Numerical continuation methods are typically used for solving the problems by setting $\rho = 1$ and solving the shooting problem. Once a solution to $\bm{\eta}^\top(t_0)$ is obtained. The value of $\rho$ is decreased, say by a factor of 0.1 or 0.5 (i.e., $\rho = \rho \times 0.1$). The previous value of $\bm{\eta}^\top(t_0)$ is used as an initial guess for this slightly-modified shooting problem. These steps are repeated until the value of $\rho$ is smaller than some prescribed user-defined value or subsequent changes in the value of cost, $J$, becomes smaller than a user-defined threshold. 
\section{Calculation of Sensitivities Using State Transition Matrix}
The TPBVPs associated with minimum-fuel trajectory optimization problems, resulting in nonlinear shooting problems that can be written as,
\begin{align}
\text{Cartesian:}~\bm{\psi}(\bm{\eta}(t_0);\rho) & = \begin{bmatrix}
\bm{r}(t_f) - \bm{r}_f \\
\bm{v}(t_f) - \bm{v}_f \\
\lambda_m(t_f)
\end{bmatrix} = \bm{0}, & \text{MEE:}~\bm{\psi}(\bm{\eta}(t_0);\rho) & = \begin{bmatrix}
\bm{x}(t_f) - \bm{x}_f \\
\lambda_m(t_f)
\end{bmatrix} = \bm{0}.
\end{align}

In these problems, the initial states, including the initial mass value, are known. Our objective is to determine the initial costate vector, denoted as $\bm{\eta}(t_0)$. While numerous off-the-shelf root-finding solvers are available, we utilize Powell's dog-leg method, a well-regarded choice for this category of problems. Powell's method is a gradient-based numerical solver that leverages the Jacobian to solve a nonlinear system of equations. If the Jacobian function is not provided, it approximates the Jacobian using a Finite Difference (FD) method. Although FD is often effective, it can introduce instability and inefficiency into the solver. Therefore, it is often advantageous to provide Powell's method with accurate values for the Jacobian of the residual error vector with respect to the unknown decision values. 

In our trajectory optimization problems, given the TPBVP, and for a fixed value of $\rho$, we intend to find the Jacobian of the final boundary condition w.r.t. the initial costates, which can be written as,
\begin{equation}
\dfrac{\partial \bm{\psi}(\bm{x}(t_f),\lambda_m(t_f),t_f)}{\partial \bm{\eta}(t_0)}.
\label{eq:boundary_to_co}
\end{equation}

To find the Jacobian, we utilize the State Transition Matrix (STM) \cite{Butcher2019}. Let $\bm{z}$ denote the states of a continuous set of differential equations denoted as, 
\begin{equation} \label{eq:zdot}
    \dot{\bm{z}}(t) = \bm{\Gamma}(\bm{z}(t),t).
\end{equation}

We intend to find the Jacobian of $\bm{z}(t_f)$ w.r.t. $\bm{z}(t_0)$, which can be written as,
\begin{equation} \label{eq:sTMdef}
    \bm{\Phi}(t_f,t_0) = \dfrac{\partial \bm{z}(t_f)}{\partial \bm{z}(t_0)}.    
\end{equation}

The solution to the differential equation that is given in Eq.~\eqref{eq:zdot} can be written as,
\begin{equation} \label{eq:abssol}
    \bm{z}(t_f) = \bm{z}(t_0) + \int_{t_0}^{t_f} \bm{\Gamma} (\bm{z}(t),t)dt.
\end{equation}

We can take the derivative of both sides of Eq.~\eqref{eq:abssol}, and using the definition given in Eq.~\eqref{eq:sTMdef}, we can write
\begin{equation}
     \bm{\Phi}(t_f,t_0) = \dfrac{\partial \bm{z}(t_f)}{\partial \bm{z}(t_0)} = I+ \int_{t_0}^{t_f} \dfrac{\partial \bm{\Gamma}(\bm{z}(t),t)}{\partial \bm{z}(t)} \dfrac{\partial \bm{z}(t)}{\partial \bm{z}(t_0)} dt = I+ \int_{t_0}^{t_f} \dfrac{\partial \bm{\Gamma}(\bm{z}(t),t)}{\partial \bm{z}(t)} \bm{\Phi}(t,t_0) dt.
\end{equation}

We proceed by taking the time derivative of both sides to obtain the differential equation that describes the time derivative of the sensitivity matrix as, 
\begin{equation}
    \dot{\bm{\Phi}}(t,t_0) = \dfrac{\partial \bm{\Gamma}(\bm{z}(t),t)}{\partial \bm{z}(t)} \bm{\Phi}(t,t_0), \quad \text{with} \quad \bm{\Phi}(t_0,t_0) = I.
    \label{eq:stm_diff_eq}
\end{equation}

In the context of the minimum-fuel trajectory optimization problems, $\bm{z}$ corresponds to the vector consisting of states $\bm{x}$, $m$, costates associated with $\bm{x}$, denoted as $\bm{\lambda}$, and the costate associated with mass, $\lambda_m$. Let $\bm{z}^\top = [\bm{x}^\top, m, \bm{\lambda}^\top, \lambda_m] \in \mathbb{R}^{14}$. From Eq.~\eqref{eq:boundary_to_co} and using the chain rule, we have
\begin{equation}
    \left [ \dfrac{\partial \bm{\psi}(t_f)}{\partial \bm{\eta}(t_0)} \right ]_{7 \times 7} = \left [\dfrac{\partial \bm{\psi}(t_f)}{\partial \bm{z}(t_f)} \right ]_{7 \times 14} \left [\dfrac{\partial \bm{z}(t_f)}{\partial \bm{z}(t_0)} \right ]_{14 \times 14} \left [\dfrac{\partial \bm{z}(t_0)}{\partial \bm{\eta}(t_0)}\right ]_{14 \times 7}.
\end{equation}

The elements of the $\partial \bm{\Gamma}(\bm{z}(t))/\partial \bm{z}(t)$ are given in Appendix A. One of the key implementation steps in this derivation is that, the state transition matrix is obtained after substituting the smooth optimal throttle into the right-hand side of the differential equations to respect the \textit{continuous} state assumption in the STM derivation. The idea of substituting the smooth control is first introduced in \cite{taheri2016enhanced} to avoid developing an event-detection logic as part of the solution of the resulting TPBVPs. 

Please note that if the smooth control is not substituted in the right-hand side of the differential equations, an event-detection logic is required to determine the precise zero-crossings of the throttle switching function to update not only the control value, but also to determine a transition matrix that is required to updating the value of the STM matrix after the control switch has occurred. The details involved in deriving the relations are given in \cite{prussing2010primer}.  

Finally, after forming the equations for both $\dot{\bm{z}}$ and $\dot{\bm{\Phi}}$, the system of differential equations consists of 210 differential equations (7 for states and mass, 7 for costates, and $14 \times 14 = 196$ for the Jacobian matrix $\bm{\Phi}$). For a randomly initialized costates, we obtain the Jacobian and boundary conditions and reiterate over initial costates using Powell's method to search for the solution that satisfies the boundary conditions. The solution procedure is similar to \cite{arya2019hyperbolic}.
\section{Results}
We investigate two benchmark fixed-time minimum-fuel trajectory optimization problems: the Earth-to-Mars and Earth-to-Dionysus problems \cite{taheri2016enhanced}. Given the bang-off-bang nature of optimal control in these scenarios, numerical solvers and integrators face significant challenges. To address this, we compare two smoothing functions, HTS and L2 (short for L2-norm-based) methods, solving both using Cartesian and MEEs as outlined in earlier sections. Additionally, we examine the impact of calculating the sensitivities using the STM for both regularization methods and for both problems. The time and distance units used to solve the problems, respectively, are $ \text{TU} = 3.1536 \times 10^7$ s , $\text{AU} = 1.496 \times 10^8$ km. The costates for the Cartesian coordinate are randomly generated in the range $[0,1]$, and for the MEEs, they are randomly generated in the range $[0,0.1]$ except for $\lambda_m$, which was randomly generated in the range $[0,1.0]$. Moreover, the relative and absolute tolerances for the numerical integrator are both set to $1.0 \times 10^{-13}$.

The code was developed using Python. The Sympy library \cite{10.7717/peerj-cs.103} is employed for modeling the problem, while SciPy \cite{2020SciPy-NMeth} is used for numerical integration and root finding. It is noteworthy that the default nonlinear root finder in SciPy utilizes the original Minpack \cite{more1980user} hybrid process, which is an implementation of a modified version of Powell's dog leg method \cite{powell1970hybrid}. In this study, we initially set the value of $\rho$ to 1 to find an initial solution. The obtained solution is then fed back into the problem with a reduced value of $\rho = 0.1\rho$ and solved again. This iterative process continues until $\rho$ reaches $1.0 \times 10^{-5}$.

\subsection{Minimum-Fuel Earth-to-Mars Problem}
In this problem, the spacecraft should leave the Earth and rendezvous with Mars over a fixed time of flight of  $348.795$ days. Only the heliocentric phase of flight is considered, and the hyperbolic excess velocity with respect to Earth (position and velocity with subscript '$E$') and Mars (position and velocity with subscript '$M$') are assumed to be zero. The parameters and boundary conditions that are used for this problem are summarized in Table \ref{tab:spacecraft_mission_values_E2M}. Note that MEEs can be obtained from the Cartesian boundary conditions \cite{junkins2019exploration}.

\begin{table}[h!]
    \centering
    \caption{Earth-to-Mars problem: spacecraft mission parameters and initial/final conditions.}
    \begin{tabular}{|c|c|}
        \hline
        \textbf{Symbol} & \textbf{Value} (unit) \\
        \hline
        $\mu$ & $132712440018 \, (\text{km}^3/\text{s}^2$) \\
        \hline
        $m_0$ & $1000 \, (\text{kg})$ \\
        \hline
        $I_{\text{sp}}$ & $2000 \, (\text{s})$ \\
        \hline
        $T_{\text{max}}$ & $0.5 \, \text{N}$ \\
        \hline
        $(\mathbf{r}_E, \mathbf{v}_E)$ & $\begin{aligned} \mathbf{r}_E &= [-140699693, -51614428, 980]^\top \, (\text{km}), \\ \mathbf{v}_E &= [9.774596, -28.07828, 4.337725 \times 10^{-4}]^\top \, (\text{km/s}) \end{aligned}$ \\
        \hline
        $(\mathbf{r}_M, \mathbf{v}_M)$ & $\begin{aligned} \mathbf{r}_M &= [-172682023, 176959469, 7948912]^\top \, (\text{km}), \\ \mathbf{v}_M &= [-16.427384, -14.860506, 9.21486 \times 10^{-2}]^\top \, (\text{km/s}) \end{aligned}$ \\
        \hline
        $\text{TOF}$ & $348.795$ (days) \\
        \hline
    \end{tabular}
    
    \label{tab:spacecraft_mission_values_E2M} 
\end{table}

\begin{figure}[htb!]
    \centering
    \includegraphics[width=0.8\linewidth]{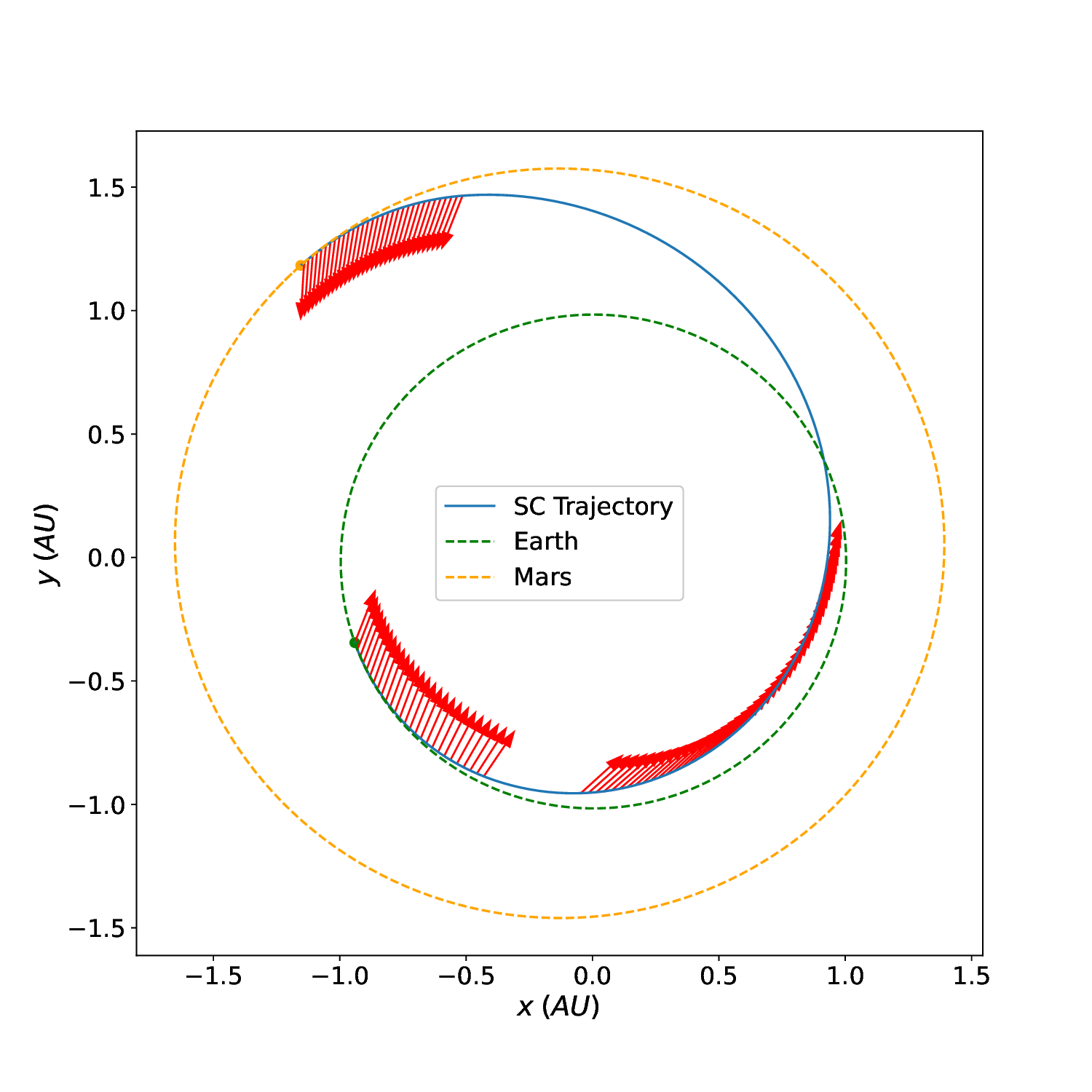}
    \caption{Earth-to-Mars problem: $x-y$ view of the optimal trajectory for $\rho = 1.0 \times 10^{-5}$. The red arrows show the thrust direction.}
    \label{fig:e2m_traj}    
\end{figure}

\begin{figure}[ht!]
    \centering
    \includegraphics[width=0.9\linewidth]{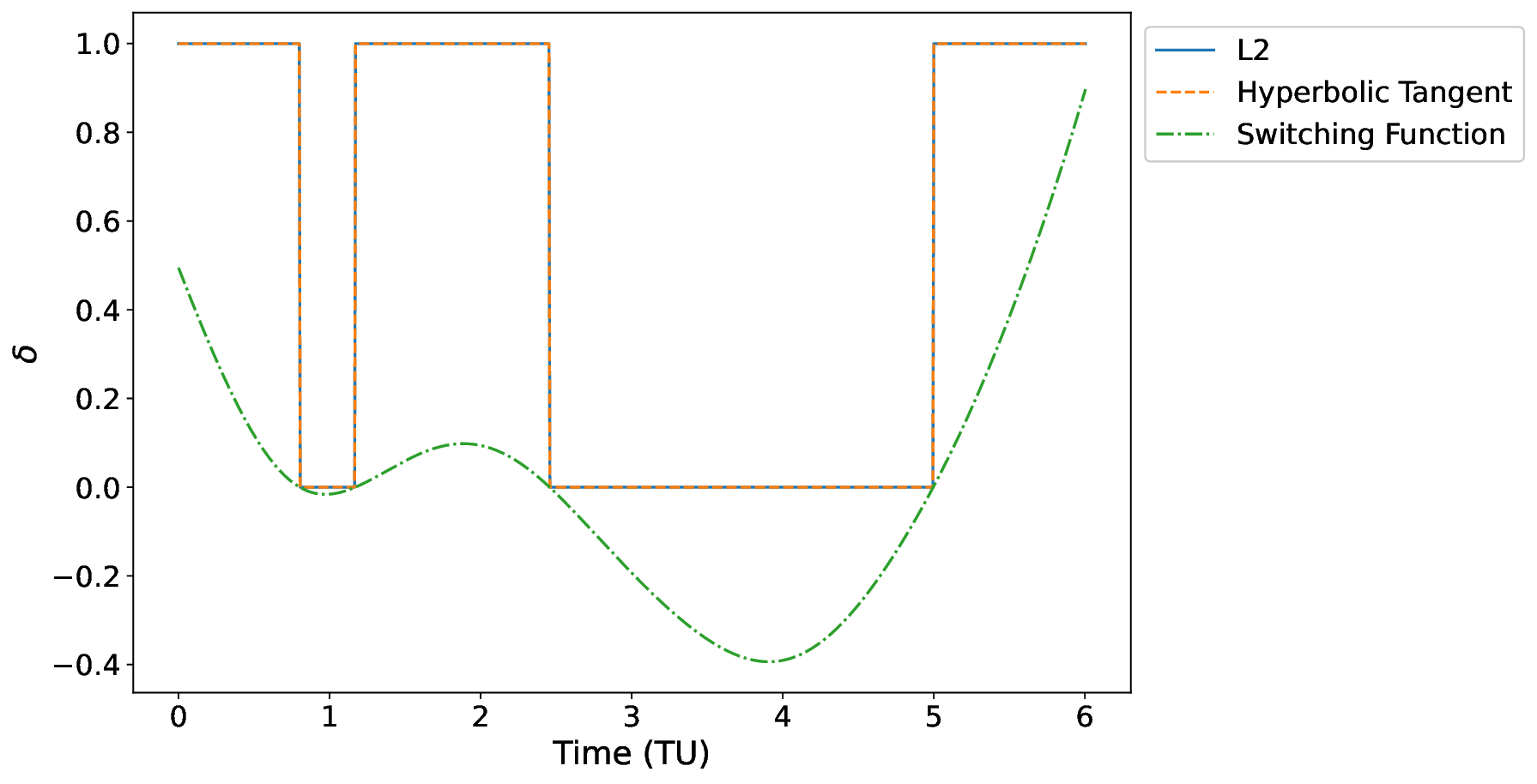}
    \caption{Earth-to-Mars problem: throttle and switching function vs. time for $\rho = 1.0 \times 10^{-5}$.}
    \label{fig:e2m_thrust}    
\end{figure}

\begin{figure}[ht!]
    \centering
    \includegraphics[width=0.9\linewidth]{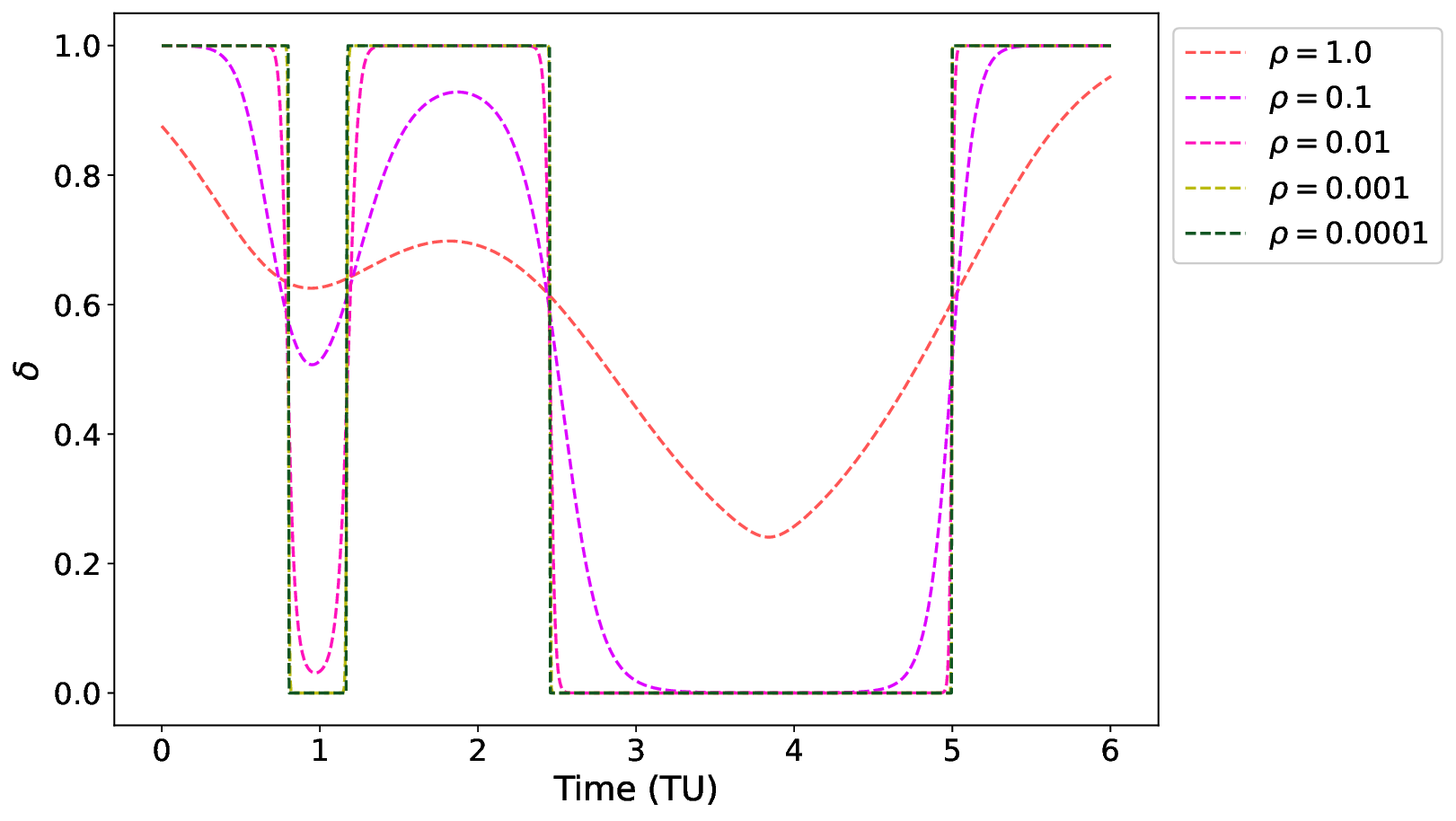}
    \caption{Earth-to-Mars problem: throttle profile vs. time for different values of $\rho$.}
    \label{fig:e2m_thrust_rho}    

\end{figure}

The result of the numerical simulations, for 100 randomly generated set of costates, are summarized in Table \ref{tab:comparison_E2M}. The last column refers to the average time it took for the solution to be found. Note that this average time includes the time for all values of \(\rho\) from 1.0 to \(10^{-5}\). From the results, it is observed that L2 Smoothing has a slight advantage over the HTS method. Moreover, confirming the results of previous studies \cite{taheri2016enhanced,arya2019hyperbolic}, using the STM method greatly improves the convergence rate of the nonlinear root-finding algorithms. In Table \ref{tab:comparison_E2M},  The optimal obtained mass for this problem is $m(t_f) = m_f = 603.935$ kg, which is consistent with the globally optimal solution for this benchmark problem \cite{taheri2016enhanced}. The minimum-fuel trajectory is shown in Fig.~\ref{fig:e2m_traj}. Figure \ref{fig:e2m_thrust} shows the throttle and switching function values for L2 and HTS methods. Figure~\ref{fig:e2m_thrust_rho} shows the throttle profile for five different values of the smoothing parameter, which is part of the numerical continuation procedure.  

\begin{table}[h]
    \centering
    \caption{Earth-to-Mars problem: comparison of convergence rate for HTS and L2 smoothing and for Cartesian and MEEs with and without the STM.}
    \begin{tabular}{|l|l|c|c|c|}
        \hline
        Smoothing Function                  & Coordinates & STM & Convergence \% & Time (s)  \\
        \hline
        Hyperbolic Tangent      & Cartesian   & True                    & 85             & 1.60                             \\
        Hyperbolic Tangent      & Cartesian   & False                   & 76             & 1.62                            \\
        \hline
        L2                      & Cartesian   & True                    & 89             & 1.60                             \\
        L2                      & Cartesian   & False                   & 78             & 1.47                             \\
        \hline
        Hyperbolic Tangent      & MEE         & True                    & 76             & 1.25                              \\
        Hyperbolic Tangent      & MEE         & False                   & 75             & 0.98                            \\
        \hline
        L2                      & MEE         & True                    & 77             & 1.37                              \\
        L2                      & MEE         & False                   & 66             & 1.29                             \\
        \hline
    \end{tabular}    
    \label{tab:comparison_E2M}
\end{table}

\begin{figure}[ht!]
    \centering
    \includegraphics[width=0.8\linewidth]{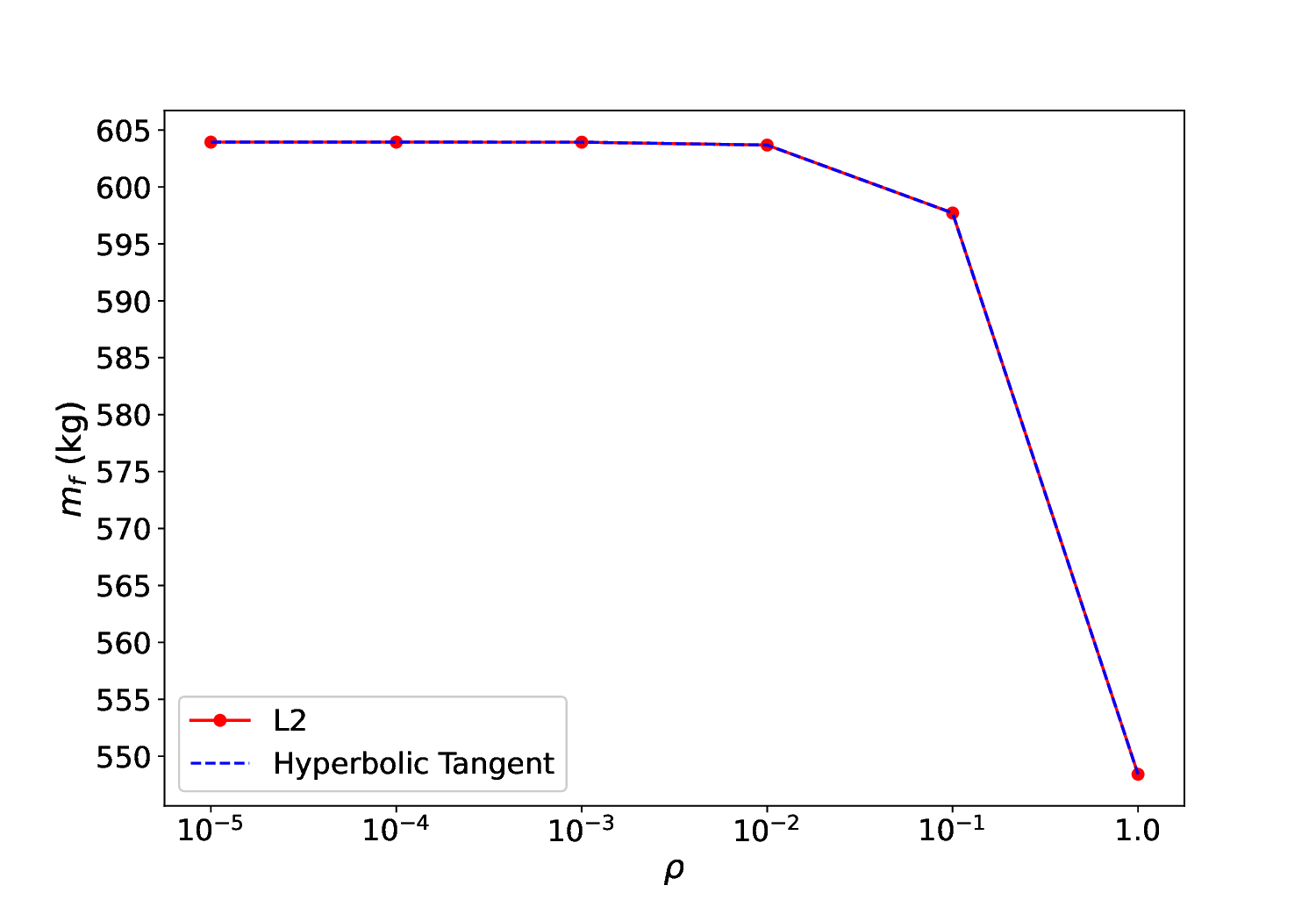}
    \caption{Earth-to-Mars problem: final mass vs. $\rho$.}
    \label{fig:e2m_mf}    
\end{figure}

\begin{figure}
    \centering
    \vspace{-100pt}

    \includegraphics[width=0.9\linewidth]{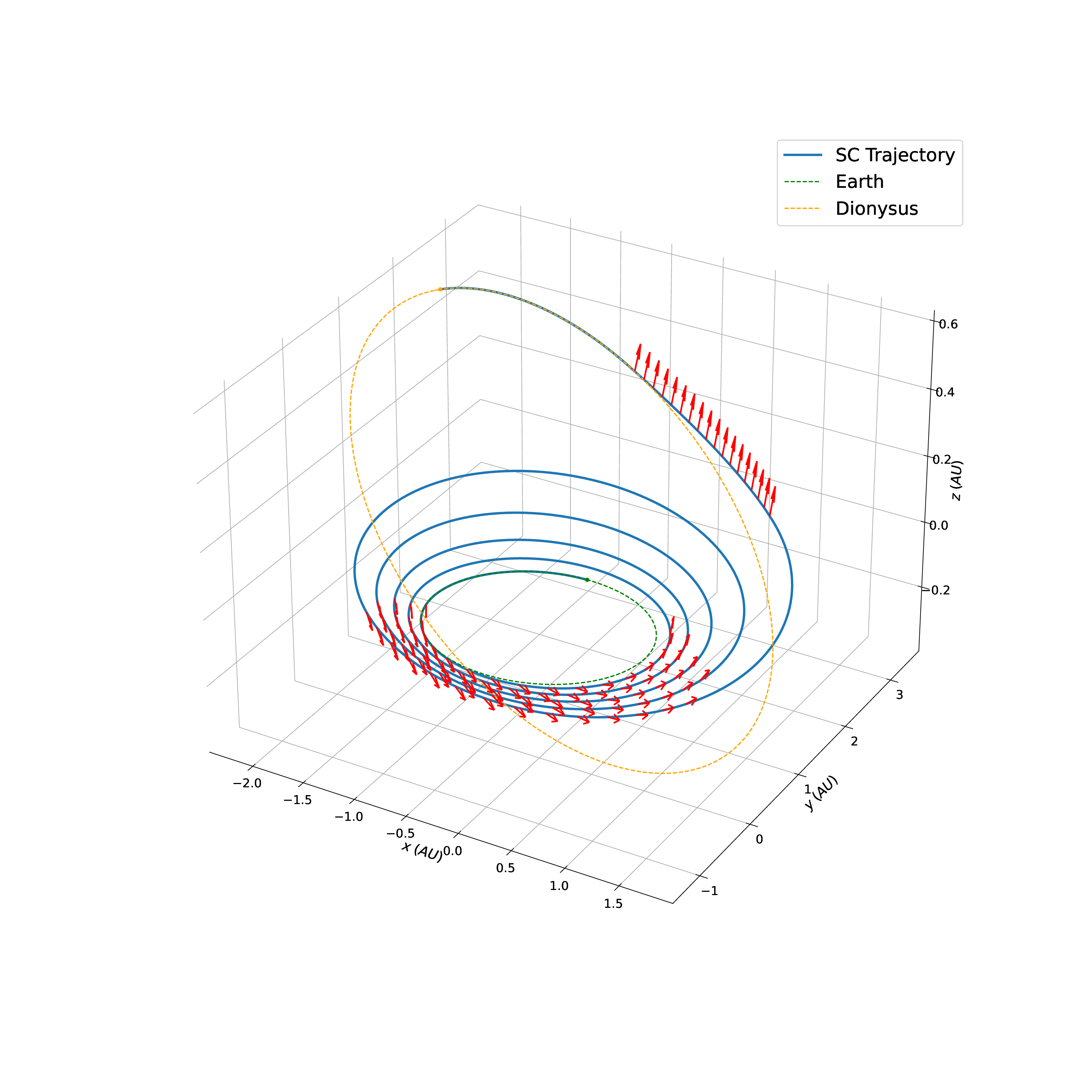}
     \vspace{-50pt}
    \caption{Earth-to-Dionysus problem: optimal trajectory for  $\rho = 1.0 \times 10^{-5}$. The red arrows show the thrust direction.}
    \label{fig:e2d_traj_3d}    
\end{figure}

\subsection{Minimum-Fuel Earth-to-Dionysus Problem} 
\begin{figure}
    \centering
    \includegraphics[width=0.8\linewidth]{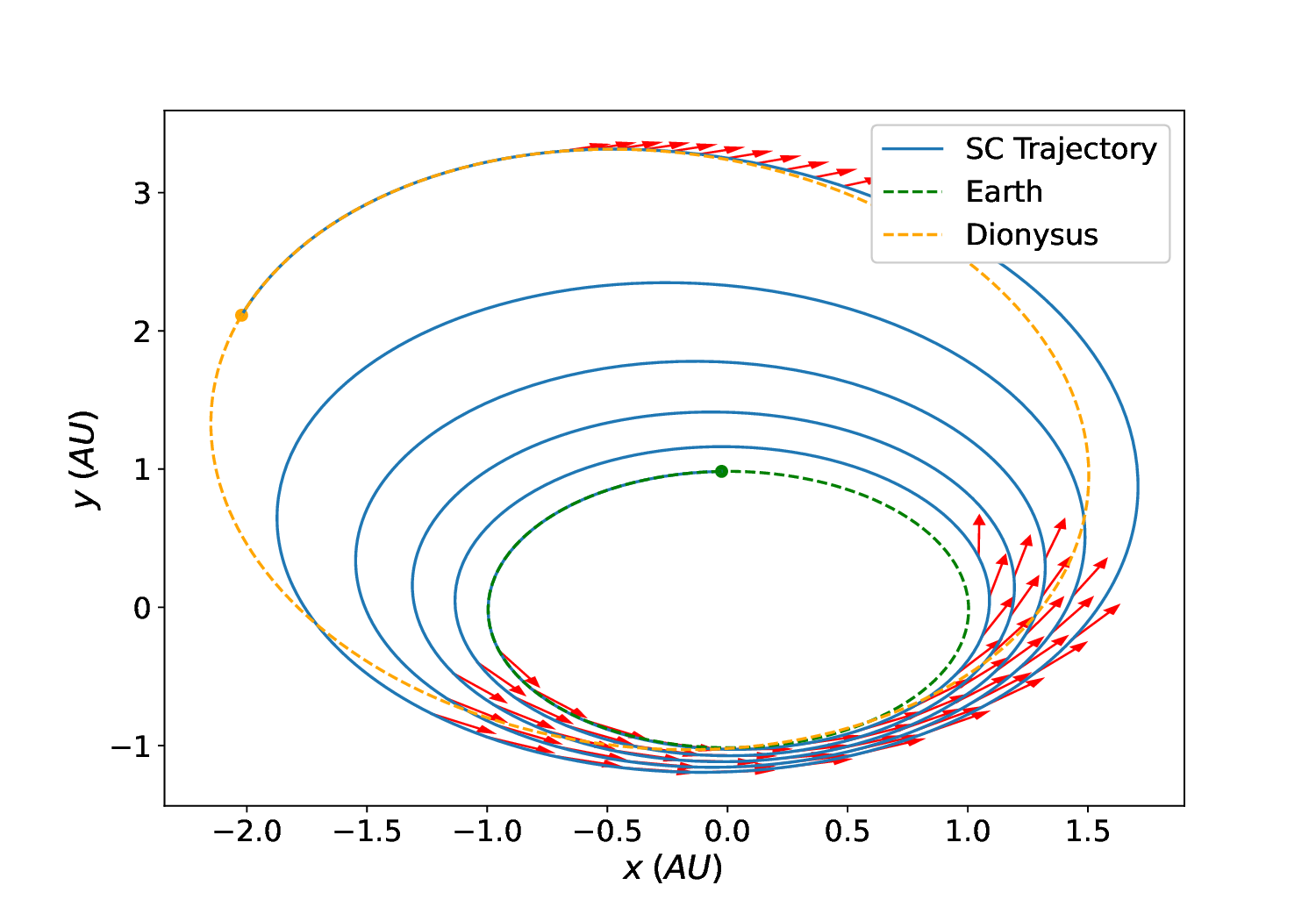}
    \caption{Earth-to-Dionysus problem: $x-y$ view of the optimal trajectory for $\rho = 1.0 \times 10^{-5}$.}
    \label{fig:e2d_traj}   
\end{figure}

\begin{figure}[ht!]
    \centering
    \includegraphics[width=0.9\linewidth]{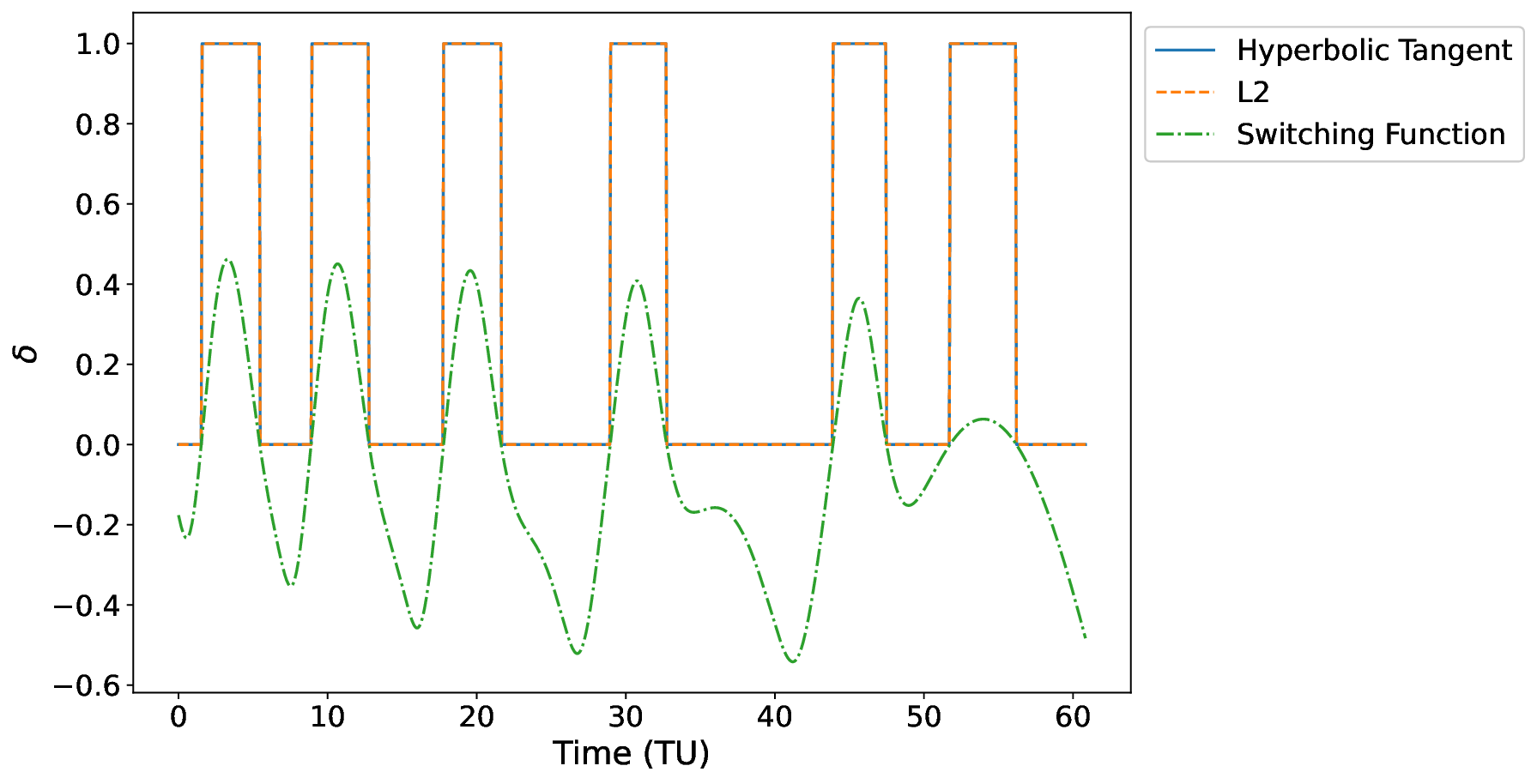}
    \caption{Earth-to-Dionysus problem: throttle and switching function vs. time for $\rho = 1.0 \times 10^{-5}$.}
    \label{fig:e2d_thrust}        
\end{figure}

\begin{figure}[ht!]
    \centering
    \includegraphics[width=0.9\linewidth]{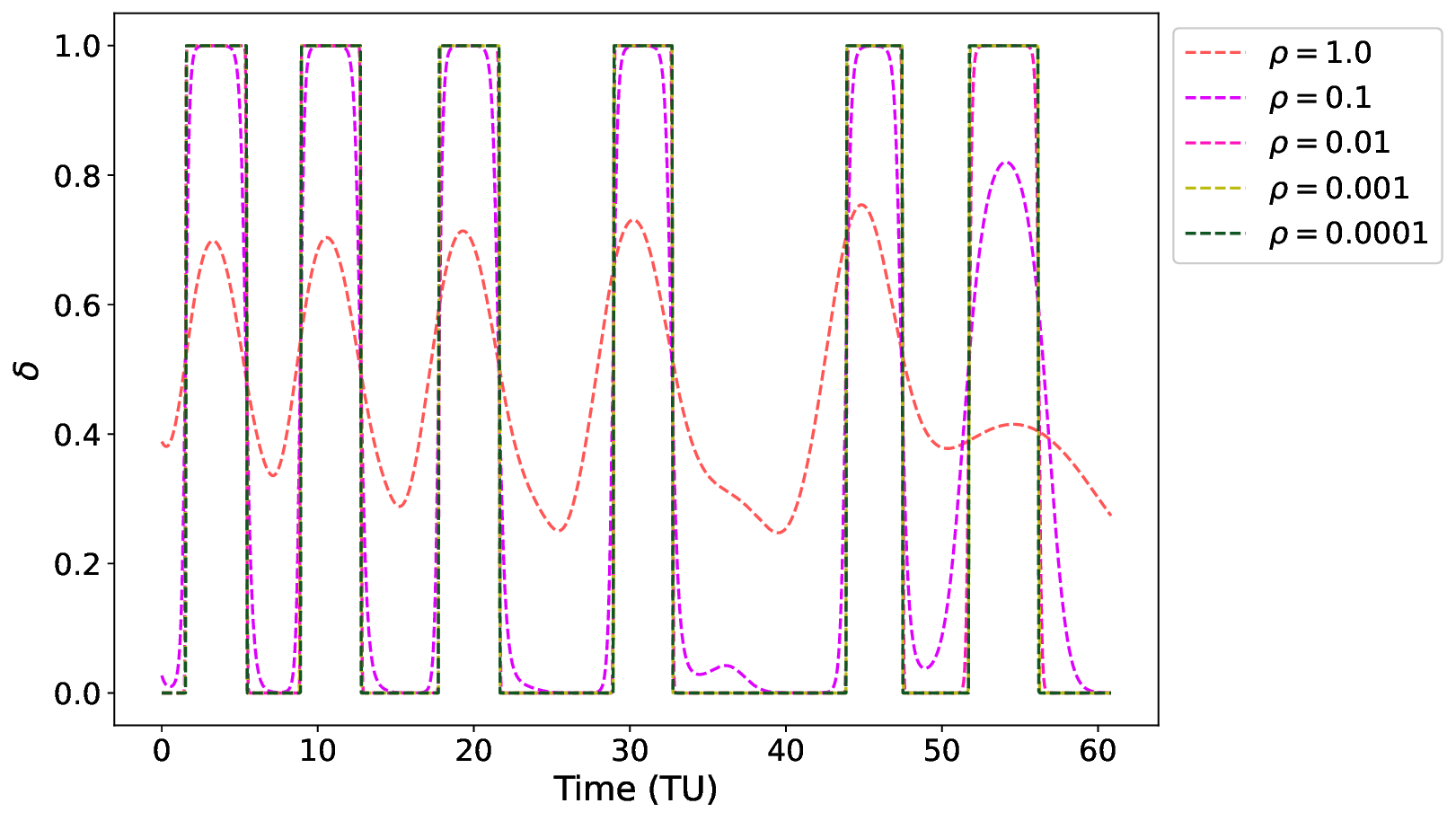}
    \caption{Earth-to-Dionysus problem: throttle vs. time for different values of $\rho$.}
    \label{fig:e2d_thrust_rho} 
\end{figure}

\begin{figure}[ht!]
    \centering
    \includegraphics[width=0.8\linewidth]{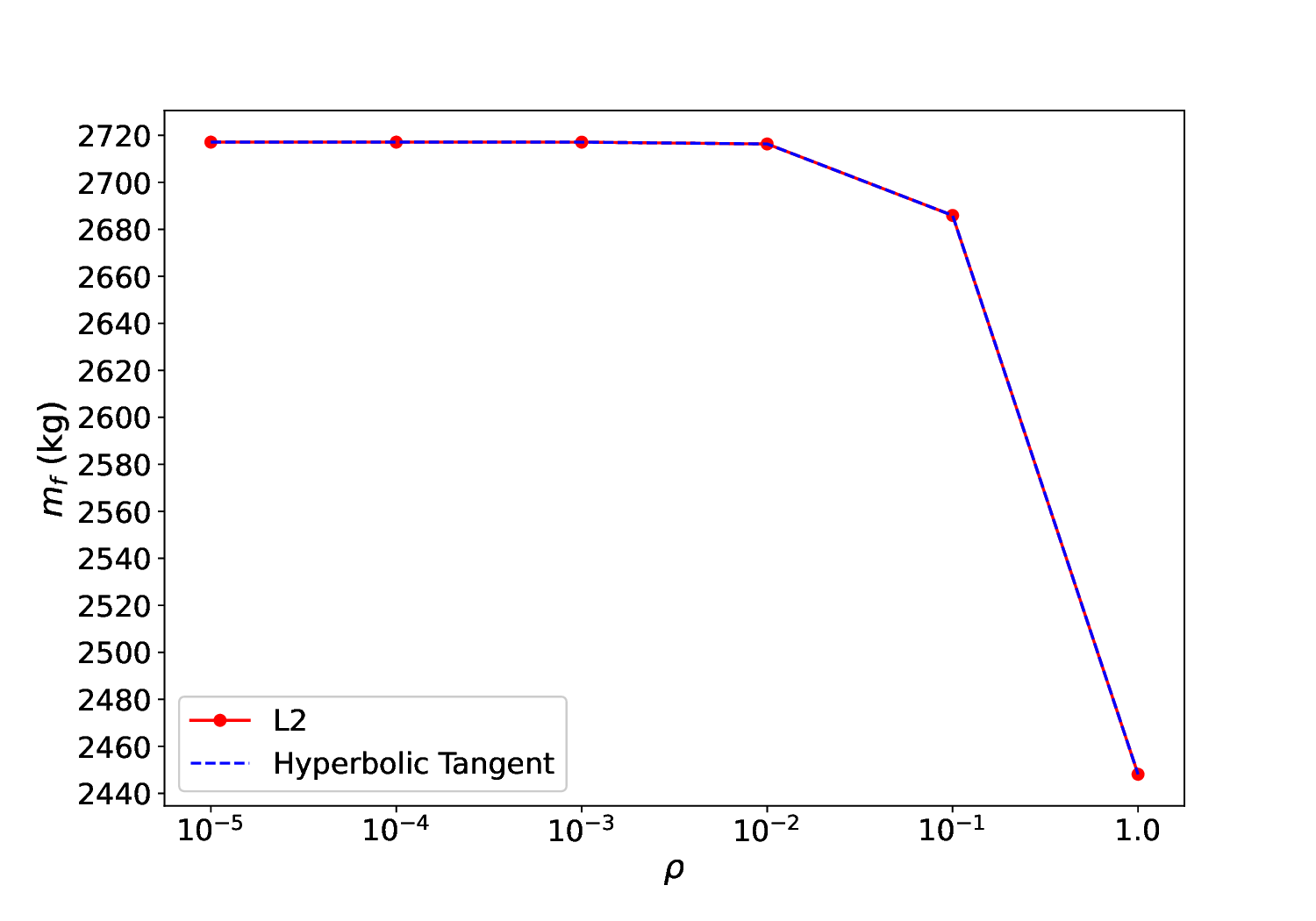}
    \caption{Earth-to-Dionysus problem: final mass vs. $\rho$.}
    \label{fig:e2d_mf} 

\end{figure}

In the Earth-to-Dionysus problem, the total flight time is 3534 days, and the parameters for this problem are summarized in Table \ref{tab:spacecraft_mission_values_E2D} with subscript ``$D$'' denoting the states for asteroid Dionysus. This problem features multiple extremal points \cite{taheri2020many} (corresponding to different number of revolutions around the Sun), which the indirect method can converge to arbitrarily when using Cartesian coordinates. The most optimal solution involves five orbital revolutions, as shown in Fig.~\ref{fig:e2d_traj_3d}. Figure \ref{fig:e2d_traj} shows the projection of the trajectory onto the $x-y$ inertial coordinate system. Figure \ref{fig:e2d_thrust} shows the time histories of the throttle and switching functions for L2 and HTS smoothing methods. Figure \ref{fig:e2d_thrust_rho} shows the throttle profile vs. time for 5 different values of the smoothing parameter.  Figure \ref{fig:e2d_mf} shows the changes in the final mass for different values of the smoothing parameter for both L2 and HTS methods. 

The results of the simulations are presented in Table \ref{tab:comparison_E2D}. In total, 100 randomly generated initial costates are considered. In this scenario, the L2 smoothing has a slight advantage over the HTS. The only exception is when the STM is not used while using MEEs. The optimal final mass for this problem is obtained as $m(t_f) = m_f = 2718.33$ kg, which is consistent with the global optimal solution reported in the literature \cite{taheri2016enhanced}.

\begin{table}[!ht]
    \centering
    \caption{Earth-to-Dionysus problem: spacecraft mission parameters and initial/final conditions.}
    \begin{tabular}{|c|c|}
        \hline
        \textbf{Symbol} & \textbf{Value} (unit) \\
        \hline
        $\mu$ & $132712440018 \, (\text{km}^3/\text{s}^2$) \\
        \hline
        $m_0$ & $4000 \, (\text{kg})$ \\
        \hline
        $I_{\text{sp}}$ & $3000 \, (\text{s})$ \\
        \hline
        $T_{\text{max}}$ & $0.32 \, (\text{N})$ \\
        \hline
        $(\mathbf{r}_E, \mathbf{v}_E)$ & $\begin{aligned} \mathbf{r}_E &= [-3637871.081, 147099798.784, -2261.441]^\top \, (\text{km}), \\ \mathbf{v}_E &= [-30.265097, -0.8486854,  0.505 \times 10^{-4}]^\top \, (\text{km/s}) \end{aligned}$ \\
        \hline
        $(\mathbf{r}_D, \mathbf{v}_D)$ & $\begin{aligned} \mathbf{r}_D &= [-302452014.884, 316097179.632, 82872290.0755]^\top \, (\text{km}), \\ \mathbf{v}_D &= [-4.533473, -13.110309,  0.656163]^\top \, (\text{km/s}) \end{aligned}$ \\
        \hline
        $\text{TOF}$ & $3534$ days \\
        \hline
    \end{tabular}
    \label{tab:spacecraft_mission_values_E2D} 
\end{table}

\begin{table}[h]
    \centering
    \caption{Earth-to-Dionysus problem: comparison of results using Cartesian and MEEs for L2 and HTS methods and with and without the STM.}
    \begin{tabular}{|l|l|c|c|c|}
        \hline
        Smoothing Function     & Coordinates & STM  & Convergence \% & Time (s)  \\
        \hline
        Hyperbolic Tangent     & Cartesian   & True  & 34             & 20.48     \\
        Hyperbolic Tangent     & Cartesian   & False & 3             & 9.33    \\
        \hline
        L2                     & Cartesian   & True  & 40             & 22.34     \\
        L2                     & Cartesian   & False & 3             & 19.87    \\
        \hline
        Hyperbolic Tangent     & MEE         & True  & 70             & 5.05     \\
        Hyperbolic Tangent     & MEE         & False & 36             & 4.23     \\
        \hline
        L2                     & MEE         & True  & 72             & 5.35     \\
        L2                     & MEE         & False & 34             & 4.01      \\
        \hline
    \end{tabular}    
    \label{tab:comparison_E2D}
\end{table}

\section{Conclusion}

This study compared the recently introduced L2-norm-based smoothing with the hyperbolic tangent smoothing method for solving minimum-fuel low-thrust trajectories using indirect optimization methods.
The resulting two-point boundary-value problems (TPBVPs) are formulated. To enhance the convergence performance of the non-linear root-finding problems associated with the resulting TPBVPs, the State Transition Matrix (STM) method is used for accurate calculation of the sensitivities. The problems were solved using Cartesian coordinates and Modified Equinoctial Orbital Elements (MEE). Two benchmark trajectory optimization problems were solved: 1) a fixed-time minimum-fuel Earth to Mars rendezvous maneuver and 2) a fixed-time minimum-fuel Earth to asteroid Dionysus rendezvous maneuver.

The results indicate that for the considered range of costates, the L2-norm-based smoothing function showed a slight advantage over the hyperbolic tangent smoothing method in terms of the percent of convergence, except in the scenario where the problem was solved using MEE without incorporating the STM method. Overall, the results demonstrated that both L2-norm-based and hyperbolic tangent smoothing methods are capable of obtaining optimal trajectories and broaden the basin of convergence of the TPBVPs associated with minimum-fuel trajectory optimization problems. 

% \section*{Appendix: B}
% None-zero elements of $\bm{SF}$ and $\bm{H}$ are
% \begin{align*}
% SF_{1,1} = & 1-\frac{1-\mu}{r^3_1} + \frac{3(1-\mu)(x+\mu)^2}{r^5_1}- \mu \frac{1}{r^3_2}+3 \mu \frac{(x+\mu-1)^2}{r^5_2},\\
% SF_{2,2} = & 1-\frac{1-\mu}{r^3_1} + \frac{3(1-\mu)y^2}{r^5_1}-\mu \frac{1}{r^3_2} +3 \mu \frac{y^2}{r^5_2},\\
% SF_{3,3} = & -\frac{1-\mu}{r^3_1} + \frac{3(1-\mu)z^2}{r^5_1}-\mu \frac{1}{r^3_2} +3 \mu \frac{z^2}{r^5_2},\\
% SF_{1,2} = & SF_{2,1} = \frac{3(1-\mu)(x+\mu)y}{r^5_1}+3 \mu \frac{(x+\mu-1)y}{r^5_2},\\
% SF_{1,3} = & SF_{3,1} = \frac{3(1-\mu)(x+\mu)z}{r^5_1}+3 \mu \frac{(x+\mu-1)z}{r^5_2},\\
% SF_{2,3} = & SF_{3,2} = \frac{3(1-\mu)yz}{r^5_1}+3 \mu \frac{yz}{r^5_2}; H_{1,2} = -H_{2,1} = 2.
% \end{align*}

%\clearpage
\section*{Appendix: A}
In this Appendix, the elements for $\dfrac{\partial \bm{\Gamma}(\bm{z}(t))}{\partial \bm{z}(t)}$ in Eq.~\eqref{eq:stm_diff_eq} are derived when Cartesian coordinates are used:

\begin{samepage}
\begin{equation}    
\frac{\partial \bm{\Gamma}}{\partial \bm{z}} = 
\begin{bmatrix}
\mathbf{0}_{3 \times 3} & \mathbf{I}_{3 \times 3} & \mathbf{0}_{3 \times 1} & \mathbf{0}_{3 \times 3} & \mathbf{0}_{3 \times 3} & \mathbf{0}_{3 \times 1} \\
\mathbf{F}_{21} & \mathbf{0}_{3 \times 3} & \mathbf{F}_{23} & \mathbf{0}_{3 \times 3} & \mathbf{F}_{25} & \mathbf{F}_{26} \\
\mathbf{0}_{1 \times 3} & \mathbf{0}_{1 \times 3} & F_{33} & \mathbf{0}_{1 \times 3} & \mathbf{F}_{35} & F_{36} \\
\mathbf{F}_{41} & \mathbf{0}_{3 \times 3} & \mathbf{0}_{3 \times 1} & \mathbf{0}_{3 \times 3} & \mathbf{F}_{45} & \mathbf{0}_{3 \times 1} \\
\mathbf{0}_{3 \times 3} & \mathbf{0}_{3 \times 3} & \mathbf{0}_{3 \times 1} & -\mathbf{I} & \mathbf{0}_{3 \times 3} & \mathbf{0}_{3 \times 1} \\
\mathbf{0}_{1 \times 3} & \mathbf{0}_{1 \times 3} & F_{63} & \mathbf{0}_{1 \times 3} & \mathbf{F}_{65} & F_{66}
\end{bmatrix}.\\
\end{equation}
\begin{align*}
&\bm{F_{21}} = \dfrac{3 \mu \bm{r} \bm{r}^\top}{\norm{\bm{r}}^5} - \dfrac{\mu}{\norm{\bm{r}}^3} \bm{I}_{3\times 3} \;,\;  \bm{F}_{23} = \dfrac{T_{\max}}{m^2 ||\bm{\lambda}_{\bm{v}}||} \bm{\lambda}_{\bm{v}} \delta(S) + \dfrac{T_{\max}c }{m^3 }\bm{\lambda}_{\bm{v}} \delta'(S) \; , \; 
\end{align*}
\begin{align*}
&\bm{F}_{25} = \dfrac{T_{\max} \bm{\lambda}_{\bm{v}} \bm{\lambda}_{\bm{v}}^\top}{m ||\bm{\lambda}_{\bm{v}}||^3 } \delta(S) - \dfrac{T_{\max} }{m ||\bm{\lambda}_{\bm{v}}||}\bm{I}_{3 \times 3} \delta(S)  -\dfrac{T_{\max}c}{(m ||\bm{\lambda}_{\bm{v}}||)^2} \bm{\lambda}_{\bm{v}} \bm{\lambda}_{\bm{v}}^\top \delta'(S) \;, \; \bm{F}_{26} = -\dfrac{T_{\max}}{m ||\bm{\lambda}_{\bm{v}}||}\bm{\lambda}_{\bm{v}} \delta'(S),  \\
&F_{33} = \dfrac{T_{\max}}{m^2}||\bm{\lambda}_{\bm{v}}|| \delta'(S) \; , \; \bm{F}_{35} = \dfrac{T_{\max}}{m ||\bm{\lambda}_{\bm{v}}||}\bm{\lambda}_{\bm{v}}^\top \delta'(S) \; , \; F_{36} = -\dfrac{T_{\max}}{c} \delta'(S),
\end{align*}
\begin{align*}
&\bm{F}_{41} = \dfrac{15 \mu \bm{r} \bm{r}^\top \bm{\lambda}_{\bm{v}} \bm{r}^\top}{||\bm{r}||^7} - \dfrac{3 \mu \bm{1}_{3 \times 1} \bm{r}^\top \bm{\lambda}_{\bm{v}}}{||\bm{r}||^5} - \dfrac{3 \mu \bm{r} \bm{1}_{1 \times 3} \bm{\lambda}_{\bm{v}}}{||\bm{r}||^5} -  \dfrac{3\mu \bm{\lambda}_{\bm{v}} \bm{r}^\top}{||\bm{r}||^5} \; , \; \bm{F}_{45} = -\dfrac{3 \mu \bm{r} \bm{r}^\top}{||\bm{r}||^5} + \dfrac{\mu}{||\bm{r}||^3} \bm{I}_{3 \times 3}\\
&F_{63} = \dfrac{2T_{\max}}{m^3} ||\bm{\lambda}_{\bm{v}}|| \delta(S) + \dfrac{T_{\max}c}{m^4}||\bm{\lambda}_{\bm{v}}||^2\delta'(S) \; , \; \bm{F}_{65} = -\dfrac{T_{\max}}{m^2 ||\bm{\lambda}_{\bm{v}}||} \bm{\lambda}_{\bm{v}}^\top \delta(S) -\dfrac{T_{\max} c }{m^3 } \delta'(S) \bm{\lambda}^{\top}, \\ &F_{66} = -\dfrac{T_{\max}}{m^2}||\bm{\lambda}_{\bm{v}}|| \delta'(S), 
\end{align*}
where $S$ is the switching function and $\delta(S)$ is the smoothing function and $\delta'(S) =  \dfrac{\partial \delta}{\partial S}$.
\end{samepage}

For each smoothing method, $\delta$ and $\delta'$ are defined as,
\begin{enumerate}
    \item Hyperbolic-Tangent-Based Smoothing:
    \begin{subequations}
    \begin{align}
        \delta_\text{tanh}(S)&=  0.5 \left[1 + \tanh\left(\frac{S}{\rho}\right) \right],\\
        \delta'_\text{tanh}(S)&=  \dfrac{0.5}{\rho}  \text{sech}^2\left(\frac{S}{\rho}\right).
    \end{align}        
    \end{subequations}
    
    \item L2-norm-Based Smoothing \cite{taheri2023l2}:
    \begin{subequations}
    \begin{align}
        \delta_\text{L2}(S) &= 0.5 \left[1 + \frac{S}{\sqrt{S^2 + \rho^2}} \right],\\
        \delta'_\text{L2}(S) &= 0.5\frac{\rho^2}{({S^2 + \rho^2})^{\frac{3}{2}}}.
    \end{align}        
    \end{subequations}
\end{enumerate}

%\newpage

\bibliographystyle{AAS_publication}
\bibliography{references}

\end{document}